\documentstyle[12pt]{article}
\begin{document}
\newtheorem{theorem}{Theorem}[section]
\newtheorem{main theorem}{Main Theorem}
\newtheorem{Satz}{Theorem}
\newtheorem{lemma}[theorem]{Lemma}
\newtheorem{lemma and definition}[theorem]{Lemma and Definition}
\newtheorem{proposition}[theorem]{Proposition}
\newtheorem{corollary}[theorem]{Corollary}
\newtheorem{definition}[theorem]{Definition}
\newtheorem{example}[theorem]{Example}
\newtheorem{examples}[theorem]{Examples}
\newtheorem{remark}[theorem]{Remark}
\newtheorem{observation}[theorem]{Observation}
\newtheorem{fact}[theorem]{Fact}
\newtheorem{facts}[theorem]{Facts}
\newtheorem{conjecture}[theorem]{Conjecture}
\newtheorem{question}[theorem]{Question}
\newcommand{\semi}{\,{\scriptstyle{\sf X}}\!{\scriptstyle{\sf I}}\,}
\newcommand{\qed}{\\ \vspace{-7mm}\begin{flushright}{\em q.e.d.}
  \end{flushright}\vspace{5mm}}
\title{On the `Section Conjecture' in anabelian geometry\thanks{I would like
to thank
Jean-Louis Colliot-Th\'{e}l\`{e}ne, Pierre D\`{e}bes, Hendrik Lenstra,
Hiroaki Nakamura, Florian Pop, Jean-Pierre Serre and Akio Tamagawa
for very helpful discussions on the subject of this paper.}}
\author{Jochen Koenigsmann\thanks{Heisenberg-Stipendiat der Deutschen
Forschungsgemeinschaft (KO 1962/1-2)}}
\date{\small {\sc Dedicated to {\bf Alexander Prestel} on the occasion of his
60th birthday}}
\maketitle\vspace{-10mm}
\begin{abstract}
Let $X$ be a smooth projective curve of genus $>1$ over a field $K$
with function field $K(X)$,
let $\pi_1(X)$ be the arithmetic fundamental group
of $X$ over $K$ and let $G_F$ denote the absolute Galois group of a field $F$.
The {\em section conjecture} in Grothendieck's anabelian geometry says that
the sections of the canonical projection
$\pi_1 (X)\rightarrow\!\!\!\!\rightarrow G_K$
are (up to conjugation) in one-to-one correspondence with the $K$-rational
points of $X$, if $K$ is finitely generated
over ${\bf Q}$. The birational variant conjectures a similar correspondence
w.r.t. the sections  of the projection
$G_{K(X)}\rightarrow\!\!\!\!\rightarrow G_K$.

So far these conjectures were a complete mystery except for the
obvious results over separably closed fields and some non-trivial
results due to Sullivan and Huisman over the reals.
The present paper proves --- via model theory ---
the birational section conjecture for all local fields
of characteristic $0$ (except ${\bf C}$),
disproves both conjectures e.g. for the fields of all real or p-adic
{\em algebraic} numbers, and gives
a purely group theoretic characterization of the sections
induced by $K$-rational points of $X$ in the birational
setting over almost arbitrary fields.

As a biproduct we obtain Galois theoretic criteria for radical solvability
of polynomial equations in more than one variable, and for a field
to be PAC, to be large, or to be Hilbertian.
\vspace{2mm}
\end{abstract}
\tableofcontents
\section{Introduction}
\subsection{The arithmetic fundamental group}
Let $X$ be a smooth (always absolutely irreducible, i.e. geometrically
connected) curve over a field $K$ with function field $K(X)$ of genus $g_X$.
For any field extension $F/K$, write $X_F:=X\otimes_K F$ for the curve
considered over $F$ and denote the set of $F$-rational points of $X$ by
$X(F)$. Fix an algebraic closure $\overline{K}$ of $K$ and the separable
closure $K^{sep}$ and the perfect hull $K^{perf}$ of $K$ in $\overline{K}$.
Let $\tilde{X}$ be the smooth completion of $X$.
Then $\tilde{X}(\overline{K})\setminus X(\overline{K})$
is a finite set, say of cardinality $n_X$.
Let $\widetilde{K(X)^X}$ denote the maximal Galois extension of $K(X)$
which is unramified over $X$.
Then the {\bf arithmetic fundamental group} $\pi_1(X)$ of $X/K$
is defined as the Galois group of $\widetilde{K(X)^X}/K(X)$.
Denoting the absolute Galois group of a field $F$ by $G_F=Gal(F^{sep}/F)$,
one obtains the following canonical exact sequences
(with commuting squares):
$$\begin{array}{ccccccccc}
1 & \longrightarrow & G_{K^{sep}(X)} & \longrightarrow &
G_{K(X)} & \stackrel{pr_{X/K}}{\longrightarrow} & G_K & \longrightarrow & 1\\
 & & \downarrow & & \downarrow & & \| & &\\
1 & \longrightarrow & \pi_1(X_{K^{sep}}) & \longrightarrow & \pi_1(X) &
\stackrel{pr^1_{X/K}}{\longrightarrow} & G_K & \longrightarrow & 1
\end{array}$$
Passing from $K$ to $K^{perf}$ leaves the diagram unchanged ($G_{K^{perf}}
\cong G_K$ etc.).
$\pi_1(X_{\overline{K}})$ ($\cong \pi_1(X_{K^{sep}})$), the arithmetic
fundamental group of $X$ over $\overline{K}$, is called {\bf geometric
fundamental group}.
If $char\, K=0$, $X$ may be regarded as curve over ${\bf C}$ and then
the geometric fundamental group is the profinite completion of the
`algebraic fundamental group' of $X_{{\bf C}}$, i.e.
the usual fundamental group of the corresponding $n_X$-fold punctured Riemann
surface which is generated by elements $\alpha_1, \beta_1,\ldots ,\alpha_{g_X}
,\beta_{g_X}, \gamma_1,\ldots ,\gamma_{n_X}$ subject to the single relation $$
\alpha_1\beta_1\alpha_1^{-1}\beta_1^{-1}\cdots \alpha_{g_X}\beta_{g_X}\alpha_{
g_X}^{-1}\beta_{g_X}^{-1}\gamma_1\cdots \gamma_{n_X}=1.$$
So if $n_X\geq 1$, then
$\pi_1 (X_{\overline{K}})$ is a free profinite group in $2g_X+n_X-1$
generators.

If $char\,K >0$ and $n_X\geq 1$, then $\pi_1(X_{\overline{K}})$ is no longer
free, but still projective ([MM], V. Theorem 5.3).
\subsection{Grothendieck's anabelian geometry}
The fundamental conjecture of Grothendieck's program called `anabelian
geometry'  says in its simplest form that a smooth hyperbolic curve $X$
over a finitely generated extension $K$ of ${\bf Q}$ is up to $K$-isomorphism
determined by its fundamental group, or, more precisely, by the projection
$pr_{X/K}^1$ (recall that $X$ is {\bf hyperbolic} if, in the notation above,
$\chi (X):=2-2g_X-n_X<0$, i.e. if $\pi_1(X)$ is non-abelian: hence the
term `anabelian'). This conjecture was proved in a series of papers by
Nakamura ([Na1]), Tamagawa ([Ta]) and Mochizuki ([Mo]), and generalised
in many ways, for example in the following
three respects: the constant field $K$ may be any {\bf sub-$p$-adic field},
i.e. any subfield of a finitely generated extension of ${\bf Q}_p$,
the fundamental group may be replaced by the quotient obtained by passing
to the maximal pro-$l$-quotient of the kernel of $pr_{X/K}^1$, and
the isomorphism version is replaced by a much more general `Hom'-version.
For an excellent survey see [MNT].

The birational version that $X$ be encoded in $G_{K(X)}$ up to $K$-birational
equivalence was proved by Pop for arbitrary curves $X$ over finitely
generated extensions of ${\bf Q}$ ([Po2]), where one should note that
the result is `absolute' in the sense that it suffices to look at $G_{K(X)}$
and get $pr_{X/K}:G_{K(X)}\rightarrow\!\!\!\!\rightarrow G_K$ for free.

In all these approaches it is clear that the Galois theoretic data
encode, in particular, $K$-rational points.
The section conjecture says that the group theoretic code for $K$-rational
points of $X$ is as simple as possible: The $K$-rational points,
according to the conjecture, should be in 1-1 correspondence with the
conjugacy classes of sections of $pr_{X/K}^1$ resp. certain canonical families
of conjugacy classes of sections of $pr_{X/K}$. We will explain this
now more precisely.
\subsection{The section property}
If $X$ is a smooth curve over a perfect field $K$, then $\widetilde{K(X)^X}$
is the intersection of all inertia subfields of $K(X)^{sep}/K(X)$ w.r.t. all
(valuations on $K(X)$ corresponding to the) points $P\in X(\overline{K})$.
Hence, the inertia subgroups of $\pi_1(X)$ w.r.t. points in $X(\overline{K})$
are trivial and so any $K$-rational point
$P\in X(K)$ canonically induces a conjugacy class $[s_P]$ of sections
$s_P$ of $pr^1_{X/K}$, where the image $s_P(G_K)$ is a decomposition subgroup
of $\pi_1(X)$ w.r.t. $P$.
If $X$ is a smooth curve over an arbitrary field $K$, this also holds for
any $P\in X(K^{perf})$ via the canonical restriction isomorphisms:
$$\begin{array}{ccc}
\pi_1 (X_{K^{perf}}) & \stackrel{pr^1_{X/K^{perf}}}{\longrightarrow} &
G_{K^{perf}}\\
\downarrow\cong & & \downarrow\cong\\
\pi_1(X_K) & \stackrel{pr^1_{X/K}}{\longrightarrow} & G_K\end{array} $$
\begin{definition}
For a smooth curve $X$ over a field $K$ we define the
{\bf section property}
$$\begin{array}{rl}
\mbox{{\bf SP}}(X/K): & \begin{array}{rccl}
\sigma_{X/K}: & X(K^{perf}) &\rightarrow &\{\mbox{non-branch sections of }
pr^1_{X/K}\}/\mbox{conj.}\\
 & P & \mapsto & [s_P]\end{array}\\
 & \mbox{is bijective,}\end{array}$$
where a section is called non-branch section, if the image is not in a
decomposition subgroup of $\pi_1(X)$ w.r.t. a branch point, i.e. a point in
$\tilde{X}\setminus X$.
\end{definition}
The section property has been conjectured by Grothendieck for complete curves
of genus $>1$ over fields which are finitely generated over ${\bf Q}$
([G], p.5).
\begin{remark}
\begin{enumerate}
\item
If $K$ is separably closed then for any smooth complete curve $X/K$,
$\sigma_{X/K}$ is obviously surjective, but not injective. Due to a result
of Sullivan ([Su]) and Huisman ([Hu]), the same holds for $K={\bf R}$:
the sections of $pr^1_{X/K}$
correspond to the connected components of $X({\bf R})$.
\item
By Theorem 19.1 of [Mo], $\sigma_{X/K}$ is injective
for any smooth complete curve $X$ of genus $>1$ over a {\bf sub-$p$-adic}
field $K$ (i.e. a subfield $K$ of a finitely generated extension
of ${\bf Q}_p$).
\item
Any branch point $P\in \tilde{X}(K^{perf})\setminus X(K^{perf})$ also induces
sections of $pr_{X/K}^1$ induced by the corresponding sections of
$pr_{X/K}$ which will be described in the next section.
\end{enumerate}
\end{remark}
\subsection{The birational section property}
For the birational version of the section property we may as well assume
that the curves to be considered are complete.
A point $P\in X(K^{perf})$ of a smooth complete curve $X$ over a field $K$
also induces sections $s$ of $pr_{X/K}$, where the image $s(G_K)$
is a complement of the inertia subgroup of a decomposition subgroup of
$G_{K(X)}$ w.r.t. $P$. Such complements always exist
([KPR]), but they need not all be conjugate. We call such sections induced
by points in $X(K^{perf})$ {\bf geometric}.
\begin{observation}
\label{geometric}
Let $X$ be a smooth complete curve over a field $K$
and let $s$ be a section of $pr_{X/K}$.
Then $s$ is geometric iff $s(G_K)$ is contained in a decomposition
subgroup of $G_{K(X)}$ w.r.t. some $P\in X(\overline{K})$.
\end{observation}
{\em Proof:}
If $s$ is geometric then, by definition, $s(G_K)$ is contained in a
decomposition subgroup of $G_{K(X)}$ w.r.t. some
$P\in X(K^{perf})\subseteq X(\overline{K})$.

For the converse, assume $s(G_K)\subseteq D_P$
for some decomposition subgroup $D_P$ of $G_{K(X)}$ w.r.t. some
(place $\overline{P}$ of $K(X)^{sep}$ above the place of $K(X)$
corresponding to) $P\in X(\overline{K})$.
Write $F$ for the fixed field of $s(G_K)=G_F$.
Since $s$ is a section of $pr_{X/K}:\, G_{K(X)}\rightarrow\!\!\!\!\rightarrow
G_K$, $K$ is relatively algebraically closed in $F$, and so,
in particular, in the fixed field $F_P$ of $D_P$.
As $P\in X(\overline{K})$ and as $F_P$ is henselian
(w.r.t. $\overline{P}$) this means $P\in X(K^{perf})$.
The fixed field of the inertia subgroup of $D_P$ is then
$K^{sep}F_P$, and $s(G_K)=G_F$ is a complement of $G_{K^{sep}F_P}$
in $D_P=G_{F_P}$:
$F\cap K^{sep}F_P=F_P$ and $F\,(K^{sep}F_P)=F_P^{sep}=F^{sep}$.$\Box$

Note that, if $K$ is not separably closed, the pendant to the injectivity of
$\sigma_{X/K}$ is always given in the birational setting: any conjugates of
two decomposition subgroups of $G_{K(X)}$ corresponding to distinct points
have trivial intersection.
\begin{definition}
For a smooth complete curve $X$ over a field $K$ we define the
{\bf birational section property}
$$\mbox{{\bf BSP}}(X/K):\mbox{ all sections of }pr_{X/K}\mbox{ are geometric}.
$$
\end{definition}
\subsection{The weak (birational) section property}
Note that the section property and the birational section property
imply that sections can only exist when there are $K^{perf}$-rational points.
We also name this weaker property:
\begin{definition}
Let $X$ be a smooth resp. a smooth complete curve over a field $K$.
Then the {\bf weak section property} resp. the {\bf weak birational
section property} are:
$$\begin{array}{lccc}
\mbox{{\bf sp}}(X/K): & \exists \mbox{ non-branch sections of }pr^1_{X/K} &
\Longleftrightarrow & X(K^{perf})\neq \emptyset\\
\mbox{{\bf bsp}}(X/K): & \exists \mbox{ sections of }pr_{X/K} &
\Longleftrightarrow & X(K^{perf})\neq \emptyset\end{array}$$
\end{definition}
\begin{remark}
\begin{enumerate}
\item
Since sections of $pr_{X/K}$ induce sections of $pr^1_{X/K}$,
{\bf sp}$(X/K)$ always implies {\bf bsp}$(X/K)$
for a smooth complete curve $X/K$.

However, in general,
for fixed $X/K$, no implications can be made between {\bf SP}$(X/K)$
and {\bf BSP}$(X/K)$, since the projection $G_{K(X)}\rightarrow\pi_1(X)$
may, in general, not split.

Yet, since $pr_{X/K}=\lim_\leftarrow pr^1_{X^\prime /K}$,
where the inverse limit is taken over all Zariski open
$X^\prime\subseteq X$, {\bf SP}$(X^\prime/K)$ for all smooth curves $X^\prime$
over $K$ implies {\bf BSP}$(X/K)$ for all smooth complete curves $X$ over $K$.
\item
{\bf bsp}$(X/K)$ holds for any smooth complete curve $X$ of genus $0$ over any
field $K$ of characteristic $\neq 2$, because any section $s$
of $pr_{X/K}$ induces an embedding of the $2$-torsion part of Brauer groups
$Br_2(K)\hookrightarrow Br_2(K(X))$ (and, more general, an embedding
of the corresponding cohomology groups $H^n(G_K)\hookrightarrow H^n(G_{K(X)})$
with coefficients in $\mu(\overline{K})$
for any $n$: $pr_{X/K}\circ s=id_{G_K}$), and because the element in $
Br_2(K)$ corresponding to a conic over $K$ is nontrivial iff $X(K)=\emptyset$.

We do not know what happens for fields of characteristic $2$.
\item
Non-existence of sections of $pr_{X/K}^1$ or of $pr_{X/K}$ is an obstruction
to the existence of points in $X(K^{perf})$. In [CS], Section 2.2,
another such obstruction is studied, the so-called {\em elementary obstruction
} saying that the canonical exact sequence of $G_K$-modules
$$1\rightarrow \overline{K}^\times\rightarrow \overline{K}(X)^\times
\rightarrow \overline{K}(X)^\times /\overline{K}^\times\rightarrow 1$$
does not split. For curves of genus $0$, or, more generally, for
Severi-Brauer varieties and for principal homogenuous spaces of tori,
this is, again, the only obstruction ([CS], Example 2.2.11). This elementary
obstruction is closely related to the abelianization of our obstruction,
i.e. to the exact sequence
$$1\rightarrow \pi_1(\overline{X})^{ab}\rightarrow \pi_1(X)^{ab}\rightarrow
G_K\rightarrow 1$$
being non-split, where $\pi_1(\overline{X})^{ab}$ is the maximal abelian
quotient of $\pi_1(\overline{X})$ and $\pi_1(X)^{ab}$ is the corresponding
extension of $G_K$ by $\pi_1(\overline{X})^{ab}$. For details cf.
[HS], Section 3.4.
\item
If $X/K$ is a smooth complete curve with $X(K^{perf})=\emptyset$
and if $G_K$ is projective then both $\mbox{\bf sp}(X/K)$ and
$\mbox{\bf bsp}(X/K)$ do not hold, as any epimorphism onto a
projective group splits. As an example you may take
$K={\bf C}(t)$ or $K={\bf C}((t))$ and consider the curve defined by
$X^3+tY^3=t^2Z^3$.
\end{enumerate}
\end{remark}
The following observation is implicit in [Ta], Prop. 2.8:
\begin{lemma}
\label{Tamagawa}
Let $K$ be a finite extension of ${\bf Q}_p$ or a field finitely generated
over {\bf Q} or a finite field. Then:
$$[\forall X/K \mbox{ {\bf (B)SP}}(X/K)]\Longleftrightarrow
[\forall X/K\mbox{ {\bf (b)sp}}(X/K)],$$
where the quantification is over all smooth (complete) curves $X/K$. For
$K={\bf R}$ the same equivalence holds, but only in the birational version.
\end{lemma}
{\em Proof:}
Let $X/K$ be a smooth complete curve, let $s$ be a section of
$pr_{X/K}^1$ resp. $pr_{X/K}$, and let $F$ be the fixed field of
$s(G_K)$.
Assuming the right hand side we have to show that $s$ comes from
a point in $X(K)$.

$F/K$ is a regular extension. So to any finite subextension $F^\prime$
of $F/K(X)$ we may choose a smooth complete curve $X^\prime$ over $K$ with
function field $F^\prime$. Denote the collection of these curves by
${\cal X}$. Then for each $X^\prime\in {\cal X}$, $s$ is also a section
of $pr_{X^\prime/K}^1$ resp. $pr_{X^\prime/K}$ and hence, by assumption,
$X^\prime (K)\neq\emptyset$. For $F^\prime\subseteq F^{\prime\prime}$
there is a canonical projection $X^{\prime\prime}\rightarrow X^\prime$
inducing a projection $X^{\prime\prime}(K)\rightarrow X^\prime (K)$.

If $K$ is a finite extension of ${\bf Q}_p$ or if $K={\bf R}$
then $X^\prime(K)$ is compact for each $X^\prime\in {\cal X}$.
And if $K$ is finite or finitely generated over ${\bf Q}$ then, by Faltings,
$X^\prime (K)$ is finite, provided $g_{X^\prime}>1$,
which holds for sufficiently large $F^\prime$. Thus,
$\lim_\leftarrow X^\prime(K)\neq\emptyset$, i.e. there is a
$K$-rational place of $F/K$.$\Box$

We expect that the Lemma is also true over any non-large field.
\section{$p$-adically closed fields}
Let us recall that a field $K$ is called {\bf $p$-adically closed} ($p$ a
prime) if $char\, K=0$ and if $K$ admits a valuation $v$
(with valuation ring ${\cal O}_v$) such that for some integer
$n>1$, $(K,v)$ is algebraically maximal with the property that
$\sharp({\cal O}_v/p{\cal O}_v)=n$. The following fact is well known
(see [PR] and [Ko1]):
\begin{fact}
\label{p-adic facts}
For a field $K$ the following conditions are equivalent:
\begin{enumerate}
\item
$K$ is $p$-adically closed.
\item
There is a (possibly trivial) henselian valuation with divisible value group
on $K$ such that the residue field is a relatively algebraically closed
subfield of some finite extension of ${\bf Q}_p$.
\item
$K$ is elementarily equivalent (in the language of fields) to some
finite extension of ${\bf Q}_p$.
\item
$G_K\cong G_F$ for some finite extension $F$ of ${\bf Q}_p$.
\end{enumerate}
\end{fact}
We shall also use the following consequence of the Galois characterization
[Ko1] of $p$-adic fields (which, in fact, already follows from
the relative characterization in [Po1]):
\begin{fact}
\label{Pop}
Let $K$ be a $p$-adically closed field, and let $F/K$ be a field extension.
If $res:\, G_F\rightarrow G_K$ is an isomorphism, then
$K$ is an elementary substructure of $F$.
\end{fact}
\subsection{The section property for $p$-adic fields}
Regarding the section property for $p$-adically closed fields we can, at the
moment, only prove the following proposition. We conjecture, however, that
the section property {\bf SP}$(X/K)$ holds for all smooth complete
curves of genus $>1$ over any local $p$-adically closed field $K$
(i.e. any finite extension $K$ of ${\bf Q}_p$).
Our proposition implies, conversely, that a $p$-adically closed field $K$
over which the section property holds for all such $X/K$ must be local:
\begin{proposition}
\label{sp}
Let $K$ be a $p$-adically closed field. Then:
\begin{enumerate}
\item
$K$ is sub-$p$-adic iff $\sigma_{X/K}$ is injective for all
smooth complete curves $X/K$ of genus $>1$.
\item
If $K$ is a proper relatively algebraically closed subfield of a
finite extension of ${\bf Q}_p$ then $\sigma_{X/K}$ is {\bf not}
surjective for any smooth complete curve $X/K$ with $X(K)\neq\emptyset$.
\end{enumerate}
\end{proposition}
{\em Proof:}
1. The direction `$\Rightarrow$' is Mochizuki's Theorem 19.1 [Mo] already
mentioned in the introduction. For the other direction, assume $K$ is not
sub-$p$-adic. Then the henselian valuation $w$ from Fact \ref{p-adic facts}.2.
is non-trivial. Let $K^{alg}=K\cap\overline{{\bf Q}}$ be the algebraic part
of $K$ and observe that $res:\, G_K\rightarrow G_{K^{alg}}$ is an isomorphism
and that $K^{alg}$ is also a relatively algebraically closed subfield
of the residue field of $w$.

Now let $X$ be a smooth complete curve
over $K^{alg}$ with $X(K^{alg})\neq\emptyset$,
say $P\in X(K^{alg})$, consider $P$ as point of $X$ over the residue field of
$w$ and lift $P$ in two different ways to points $P_1\neq P_2\in X(K)$
(via the place corresponding to $w$).
Then $P_1$ and $P_2$ induce the same section of $pr^1_{X/K}$:
$$\begin{array}{ccc}
\pi_1(X_K) & \stackrel{pr^1_{X/K}}{\longrightarrow} & G_K\\
\downarrow\cong & & res\downarrow\cong\\
\pi_1(X_{K^{alg}}) &\stackrel{pr^1_{X/K^{alg}}}{\longrightarrow} & G_{K^{alg}}
\end{array}$$

2. If $K$ is a proper relatively algebraically closed subfield of
a finite extension $F$ of ${\bf Q}_p$, and if $X$ is a curve over $K$
with $X(K)\neq\emptyset$, then there are points in $X(F)\setminus X(K)$. By
Mochizuki's injectivity result, the sections of $pr^1_{X/K}$`$=$'$pr^1_{F/K}$
induced by such points do not come from $K$-rational points.$\Box$
\subsection{The birational section property for $p$-adic and for real closed
fields}
For the birational section property we can give the complete picture
over $p$-adically closed (and real closed) fields:
\begin{proposition}
\label{bsp}
Let $K$ be a $p$-adically closed or real closed field. Then:\\
\begin{enumerate}
\item
{\bf bsp}$(X/K)$ holds for all smooth complete curves $X/K$.
\item
{\bf BSP}$(X/K)$ holds for all smooth complete curves $X/K$
iff $K$ is a local field, i.e. $K$ is a finite extension of ${\bf Q}_p$
or $K={\bf R}$.
\end{enumerate}
\end{proposition}
{\em Proof:}
1. If $X$ is a smooth complete curve over $K$ and $s$ is a section
of $pr_{X/K}:\, G_{K(X)}\rightarrow G_K$, then $s$ is an isomorphism of $G_K$
onto $G_F$ for some algebraic extension  $F$ of $K(X)$ in $\overline{K(X)}$.
Hence, by Fact \ref{Pop}, $F$ is an elementary extension of $K$
($s^{-1}=res:\, G_F\rightarrow G_K$). Since the $K$-curve $X$ has an
$F$-rational point ($K(X)\subseteq F$), it, therefore, also has a $K$-rational
point.

2. `$\Leftarrow$' follows immediately from 1. and Lemma \ref{Tamagawa}.
For `$\Rightarrow$', assume $K$ is not a local field and choose a smooth
complete curve $X$ over
$K^{alg}$ with $X(K^{alg})\neq\emptyset$.

We distinguish two cases.

{\bf Case 1:} {\em $K$ is a proper
relatively algebraically closed subfield of a local field.} Then one
can prolong the $p$-adic valuation resp. the ordering on $K$ to $K(X)$
in such a way that it remains a rank-1-valuation resp. an archimedean
ordering. Thus, the corresponding $p$-adic resp. real closure $F$ of $K(X)$
in $\overline{K(X)}$ induces a non-geometric section of
$pr_{X/K}$: $F$ cannot have a non-trivial henselian valuation
which is trivial on $K$.

{\bf Case 2:} {\em If $K$ is not a subfield of a local field}
then $K$ admits a henselian valuation $w$ with non-trivial divisible value
group $\Gamma_w$ and residue field of characteristic $0$
(cf. Fact \ref{p-adic facts}.2).
Since $X(K^{alg})
\neq\emptyset$, $w$ can be prolonged to a valuation $u$ of $K(X)$
in such a way that $\Gamma_w$ is cofinal in $\Gamma_u$.
The fixed field $F$ of a complement of the inertia subgroup of the
decomposition subgroup of $G_{K(X)}$ w.r.t. $u$ then induces a non-geometric
section of $pr_{X/K}$, for the same reason as in case 1.$\Box$
\begin{remark}
It is clear from the proof of the weak birational section property
{\bf bsp}$(X/K)$ for $p$-adically or real closed $K$ that this generalizes to
smooth complete varieties $X$ over $K$ of arbitrary dimension.

However, the strong (birational) section property has no analogue
in higher dimension: while for $\dim X=1$ any valuation on $K(X)$
which is trivial on $K$ is geometric, there may be many non-geometric
valuations on $K(X)$ which are trivial on $K$ and have residue field $K$
(e.g. with archimedean value group of rational rank $=\dim X$) if $\dim X>1$.
Such valuations induce non-geometric sections of the analogous $pr_{X/K}$.
\end{remark}
\subsection{Applications to number fields}
\subsubsection{Global sections give local points}
The following corollary is immediate from Proposition \ref{bsp}:
\begin{corollary}
Let $X$ be a smooth complete curve over a number field $K$
and assume that $pr_{X/K}$ has a section.
Then $X(\hat{K})\neq\emptyset$ for all completions $\hat{K}$ of $K$.
\end{corollary}
{\em Proof:}
Any section of $pr_{X/K}$ induces a section of $pr_{X/\hat{K}^{alg}}$.
Hence, by Proposition \ref{bsp}.1, $X(\hat{K}^{alg})\neq\emptyset$,
and so $X(\hat{K})\neq\emptyset$.$\Box$

Recall that a field is called {\bf pseudo $p$-adically or pseudo real closed}
if it satisfies a local-global-principle for rational points on varieties
w.r.t. all $p$-adic resp. all real closures.
\begin{corollary}
If $K$ is a pseudo $p$-adically or a pseudo real closed field
then {\bf bsp}$(X/K)$ holds for any curve $X$ over $K$.
\end{corollary}
\subsubsection{The local-global-principle for mere covers}
In [D], D\`{e}bes proved a local-global-principle for Galois-covers
of curves defined over number fields. For so-called mere covers
such a principle is expected to fail in general, yet no counterexamples
could be found so far. We will show that the birational section property
for all curves over number fields implies the existence of such
counterexamples.
\begin{definition}
Let $X$ be a smooth complete curve defined over a number field $K$.
Then the {\bf local-global-principle for mere covers} of $X$
is the following property:\\
{{\bf LGP}}$(X/K)$:
{\em any $\overline{K}$-cover $Y\rightarrow\!\!\!\!\rightarrow X$
which is definable over all local closures $\hat{K}$ of $K$
(i.e. there is a $\hat{K}$-cover $Y_{\hat{K}}\rightarrow\!\!\!\!\rightarrow X$
with $Y=Y_{\hat{K}}\otimes \overline{K}$)
is locally compatibly definable over $K$ (i.e. there is a $K$-cover
$Y_K\rightarrow\!\!\!\!\rightarrow X$ with $Y_{\hat{K}}=Y_K\otimes \hat{K}$
for any local closure $\hat{K}$ of $K$).}\\
Here `local closure' means a henselisation or a real closure of $K$.
\end{definition}
\begin{proposition}
Let $K$ be a number field. If the birational section conjecture
{\bf BSP}$(X/K)$ holds for any smooth complete curve $X/K$,
then there is a counterexample to the local-global-principle for mere
covers over $K$.
\end{proposition}
{\em Proof:}
We prove the contrapostion of the implication in the proposition,
and assume {\bf LGP}$(X/K)$ holds for all smooth complete curves $X$
over $K$. We have to find some $X_0/K$ with a non-geometric section
for $pr_{X_0/K}$. The trick is to choose $X_0$ to be a smooth complete
curve over $K$ with $X_0(K)=\emptyset$, but with
$X_0(\hat{K})\neq\emptyset$ for all local closures $\hat{K}$ of $K$.
We shall construct a (necessarily non-geometric) section of
$pr_{X_0/K}: G_{K(X_0)}\rightarrow G_K$.

Let $$\overline{K}(X_0)=L_0\subseteq L_1\subseteq L_2\subseteq \ldots
\subseteq \overline{K(X_0)}$$
be a tower of finite field extensions of $\overline{K}(X_0)$
such that each $L_i$ is Galois over $K(X_0)$ and such that
$\overline{K(X_0)}=\bigcup_i L_i$ (this is possible since
$\overline{K}(X_0)/K(X_0)$ is a Galois extension and so any
extension of $\overline{K}(X_0)$ has only finitely many conjugates over
$K(X_0)$).

We shall construct a sequence of $K$-covers
$$ X_0\leftarrow\!\!\!\!\leftarrow X_1\leftarrow\!\!\!\!\leftarrow X_2
\leftarrow\!\!\!\!\leftarrow\ldots$$
of curves (or, equivalently, of function fields
$K(X_0)\subseteq K(X_1)\subseteq K(X_2)\subseteq\ldots$ over $K$)
such that, for each $i$, $L_i=\overline{K}(X_i)$ and
$X_i(\hat{K})\neq\emptyset$ for any local closure $\hat{K}$ of $K$.
If this is achieved, we are done:
take $F:=\bigcup_{i=0}^\infty K(X_i)$ and observe that
$res:\, G_F\rightarrow G_K$ is an isomorphism
(surjective, since $F/K$ is regular, and injective, since
$F\overline{K}=\bigcup_{i=0}^\infty L_i =\overline{K(X_0)} =\overline{F}$);
hence $res^{-1}$ is a section of $pr_{X_0/K}$.

For the construction of the $X_i$, we start with the given $X_0$,
and show how to obtain, for $i\geq 0$, $X_{i+1}$ from $X_i$.
Since $L_{i+1}/K(X_0)$ is Galois, so is $L_{i+1}/\hat{K}(X_i)$
for any local closure $\hat{K}$ of $K$.
As $X_i(\hat{K})\neq\emptyset$, and as $\hat{K}$ is large,
we may choose $\hat{P}_i\in X_i(\hat{K})$ such that $\hat{P}_i$
is unramified in $L_{i+1}$ and let $\hat{E}_{i+1}$ be a decompostion
subfield of $L_{i+1}/\hat{K}(X_i)$ w.r.t. $\hat{P}_i$.
Then $\hat{E}_{i+1}/\hat{K}$ is a function field with $\hat{K}$-rational point
and $L_i\hat{E}_{i+1}=L_{i+1}$, i.e. the $\overline{K}$-cover corresponding
to the inclusion of $\overline{K}$-function fields $L_i\subseteq L_{i+1}$
is definable over $\hat{K}$.
Now apply {\bf LGP}$(X_i/K)$ to obtain a $K$-cover $X_{i+1}
\rightarrow\!\!\!\!\rightarrow X_i$ with $\hat{K}(X_{i+1})=\hat{E}_{i+1}$
for each local closure $\hat{K}$ of $K$, and so, in particular, with
$X_{i+1}(\hat{K})\neq\emptyset$ and $\overline{K}(X_{i+1})=L_{i+1}$.$\Box$
\section{Large countable fields}
Recall that a field $K$ is called {\bf large} (or {\bf ample})
if any variety of dimension $\geq 1$ (or, equivalently,
any curve) defined over $K$ with one smooth $K$-rational point
has infinitely many $K$-rational points.
\begin{proposition}
Let $K$ be a large countable field (e.g. $K={\bf Q}_p\cap\bar{\bf Q}$
or $K={\bf R}\cap\bar{\bf Q}$)
and let $X$ be a smooth complete curve over $K$ with $X(K)\neq\emptyset$.
Then $X/K$ has neither the section property {\bf SP}$(X/K)$
nor the birational section property {\bf BSP}$(X/K)$.
\end{proposition}
{\em Proof:}
We shall simultaneously prove the following two claims:
\begin{enumerate}
\item
If $\sigma_{X/K}$ is injective then $pr^1_{X/K}$ has a non-geometric section.
\\resp.
\item
$pr_{X/K}$ has a non-geometric section.
\end{enumerate}
It is clear that this proves the proposition.

Replacing $K$ by $K^{perf}$ we may assume that $K$ is perfect.
Let $\overline{K}(X)=L_0\subseteq L_1\subseteq \ldots$
be a tower of finite separable extensions of $\overline{K}(X)$ in
$\widetilde{K(X)^X}$ resp. in $K(X)^{sep}$ such that each $L_i$ is Galois
over $K(X)$ and such that $\bigcup_{i=0}^\infty L_i=\widetilde{K(X)^X}$ resp.
$K(X)^{sep}$. Let $X(K)=\{P_1, P_2,\ldots\}$.

We shall construct a tower of function fields
$$K(X)=K(X_0)\subseteq K(X_1)\subseteq \ldots $$ of smooth complete
curves $X_i$ over $K$ in $\widetilde{K(X)^X}$ resp. in $K(X)^{sep}$
such that, for each $i>0$,
$$\begin{array}{rl}
\star_i & L_i\subseteq \overline{K}(X_i)\\
\star\star_i & X_i(K)\neq\emptyset\\
\star\star\star_i & X_i(K)\mbox{ contains no points above }P_i
\end{array}$$
If this is achieved, then, as in the proof of the previous proposition,
$F:=\bigcup_{i=0}^\infty K(X_i)$ does the job:
$F/K$ is regular, $F\overline{K}=\bigcup_{i=0}^\infty\overline{K}(X_i)
=\bigcup_{i=0}^\infty L_i=\widetilde{K(X)^X}$ resp. $K(X)^{sep}$,
so $res:\, Gal (\widetilde{K(X)^X}/F)\rightarrow G_K$ resp.
$res:\, G_F\rightarrow G_K$ is an isomorphism,
and $res^{-1}$ is the desired section: it is non-geometric because no
$P\in X(K)$ has a $K$-rational prolongation to $F$.

For the construction, we start with $X=X_0$ and assume, for $i\geq 0$,
that $K(X_0)\subseteq K(X_1)\subseteq\ldots\subseteq K(X_i)$
have been constructed according to $\star_j -\star\star\star_j$
for each $j\leq i$.
Now choose $P\in X_i(K)$ not above $P_{i+1}$ and choose $k>i$
such that $P_{i+1}$ has no $K$-rational prolongation to a
decomposition subfield $K(X_{i+1})$ of the Galois extension
$L_kK(X_i)/K(X_i)$ w.r.t. $P$
(here the injectivity assumption enters), and such that
$P$ is unramified in $L_{i+1}/L_i$.
Then $X_{i+1}$ satisfies $\star_{i+1}-\star\star\star_{i+1}$.$\Box$
\section{Galois characterization of rational points
over almost arbitrary fields}
\subsection{Group theoretic description of geometric sections}
We apply our characterization of decomposition subgroups of absolute
Galois groups in [Ko2] to provide the `local theory' for (1-dimensional)
birational anabelian geometry over almost arbitrary fields (including
all sub-$p$-adic fields, but also all finitely generated fields of positive
characteristic, ${\bf Q}^{ab}$, ${\bf Q}^{solv}$, ${\bf Q}_p^{ab}$,
${\bf C}(t)$ etc.):
\begin{theorem}  \label{arbitrary}
Let $K$ be a field such that
\begin{itemize}
\item
$K$ is not separably closed or real closed
\item
if $char\, K=p>0$ then $G_K$ is {\em not} a pro-$p$ group
\item
$K$ either admits no non-trivial henselian valuation or
$K$ admits a henselian rank-1-valuation of mixed characteristic $(0,p)$
and $p\mid\sharp G_K$.
\end{itemize}
Let $X$ be a smooth complete curve over $K$ and let $s$ be
a section of $pr_{X/K}:\, G_{K(X)}\rightarrow\!\!\!\!\rightarrow G_K$.

Then $s$ is geometric iff $s(G_K)$ normalizes some pro-cyclic subgroup
$C$ of $G_{K(X)}$ in $G_{K(X)}$ (i.e. $\langle C,s(G_K)\rangle
=C\semi s(G_K)\leq G_{K(X)}$) with
$$C\cong\left\{\begin{array}{cl}
\hat{{\bf Z}} & \mbox{if }char\, K=0\\
\hat{{\bf Z}}/{\bf Z}_p=\prod_{q\neq p}{\bf Z}_q & \mbox{if }char\, K=p>0.
\end{array}\right.$$
\end{theorem}
{\em Proof:}
We denote by $\rho:\, G_{K^{perf}(X)}\rightarrow G_{K(X)}$ the canonical
restriction isomorphism.
If $s$ is geometric, say induced by $P\in X(K^{perf})$, then $\rho^{-1}(
s(G_K))$ is a complement of the inertia subgroup $I$ of some decomposition
subgroup $D$ of $G_{K^{perf}(X)}$ w.r.t. $P$. Now choose a complement
$C^\prime$ of the ramification subgroup of $I$ ([KPR]) and let $C=\rho(
C^\prime )$. Then $C$ is pro-cyclic
of the indicated shape and is normalized by $s(G_K)$ in $G_{K(X)}$.

For the converse, assume $C\leq G_{K(X)}$ is a pro-cyclic subgroup of the
indicated shape and normalized by $s(G_K)$ in $G_{K(X)}$.
Then, by Theorem 1 of [Ko2], the fixed field $F$ of the
subgroup $G_F=C\semi s(G_K)$ of $G_{K(X)}$
carries a tamely branching henselian valuation $v$: more precisely, for any
prime $p$ with $p^2\mid (\sharp C,\sharp G_K)$ there is a henselian valuation
$v$ {\bf tamely branching at $p$}, i.e. with non-$p$-divisible value group
(so $v$ is non-trivial) and residual characteristic $\neq p$.
Since $K$ is relatively algebraically closed in $F$, the restriction of $v$ to
$K$ is henselian. So, if $K$ has no non-trivial henselian valuation, then $v$
is trivial on $K$, i.e. it comes from a geometric place on $K(X)$, and so $s$
is by Obs. \ref{geometric} geometric. If $K$ has a henselian rank-1-valuation
of mixed characteristic $(0,p)$ and $p\mid\sharp G_K$, then $v$ may
be chosen to be tamely branching at $p$, and so the residual characteristic
of $(F,v)$ is different from $p$. This, again, forces $v$ to be
trivial on $K$, since a henselian rank-1-valuation of mixed characteristic
$(0,p)$ allows
no non-trivial henselian valuations of residual characteristic $\neq p$
on the same field. So we can proceed as above.$\Box$
\begin{remark}
It is clear that, if a pro-cyclic subgroup $C$ of $G_{K(X)}$ satisfies
the condition in the theorem, then $C$ {\em is} the inertia subgroup
of $G_F$, i.e. $C$ is contained in a (complement of the ramification
subgroup of the) inertia subgroup of some geometric decomposition subgroup
of $G_{K(X)}$: the fixed field of $s(G_K)$ is a purely tamely ramified
extension of $F$ and $K(X)^{sep}/F$ is purely inert w.r.t. $v$.
In particular, $C\leq G_{K^{sep}(X)}$.

If $X$ is a smooth curve over a field $K$ {\em which is finitely generated
over }${\bf Q}$, then Nakamura proves a similar sufficient inertia condition
for subgroups of $\pi_1(X)$ (Theorem 3.4 of [Na1]):
If $C$ is a non-trivial pro-cyclic subgroup of $\pi_1(X_{\overline{K}})$ such
that there exists a section $s$ of $pr^1_{X/K}$ with $s(G_K)$ normalizing $C$
and acting on it via the cyclotomic character of $G_K$,
then $C$ is contained in an inertia subgroup of $\pi_1(X)$ w.r.t.
a $K$-rational point in $\tilde{X}\setminus X$.
\end{remark}
Define a profinite group $G$ to be a {\bf hensel group} if there is some
Sylow subgroup $P$ of $G$ containing some non-trivial normal abelian subgroup
$N\lhd P$, but no procyclic open subgroup. Theorem 1 of [Ko2] says that
a field whose absolute Galois group is a hensel group, admits a tamely
branching henselian valuation. This allows us to give a Galois
characterization of rational points:
\begin{corollary} \label{points} Let $K$
and $X$ be as in Theorem \ref{arbitrary}. Then the map
$$\begin{array}{rcl}
X(K^{perf}) & \rightarrow & \left\{ \begin{array}{l}
\mbox{conjugacy classes of maximal hensel subgroups}\\
D\mbox{ of }G_{K(X)}\mbox{ with }pr_{X/K}(D)=G_K\end{array}\right\}\\
P & \mapsto & [D_P]\end{array}$$
is bijective.
Here $D_P$ denotes a decomposition subgroup of $G_{K(X)}$ w.r.t. $P$.
\end{corollary}
{\em Proof:} By the proof of Theorem \ref{arbitrary}, any hensel subgroup
$D$ of $G_{K(X)}$ with $pr_{X/K}(D)=G_K$
is the absolute Galois group of a henselian algebraic extension $(F,v)$ of
$(K(X),v_P)$ for some $P\in X(K^{perf})$ (with corresponding valuation $v_P$):
note that $pr_{X/K}(D)=G_K$ implies that $K$ is algebraically closed in $F$.
So $D\subseteq D_P$ for some decomposition subgroup $D_P$ of $G_{K(X)}$ w.r.t.
$P$. Moreover, if for some $P^\prime\in X(K^{perf})$, $D_{P^\prime}\subseteq D
$ then $D_{P^\prime}\subseteq D_P$, and so, by a well-known theorem of
F.K.Schmidt, $P^\prime=P$ and $D_{P^\prime}=D=D_P$. So the decomposition
subgroups of $G_{K(X)}$ w.r.t. $K^{perf}$-rational points are indeed the
maximal hensel subgroups $D$ of $G_{K(X)}$ with $pr_{X/K}(D)=G_K$,
and distinct points have non-conjugate decomposition subgroups.$\Box$

Note that if $K$ satisfies the conditions of Theorem \ref{arbitrary}, then so
does any finite extension $L/K$. So we get a Galois characterization of
$L^{perf}$-rational points for any finite extension $L/K$ and, thus, a
Galois characterization of all $P\in X(\overline{K})$.
\subsection{Back to the roots: solving equations by radicals}
If $L/K$ is a (possibly infinite) Galois extension, if $L$ (and hence $K$)
satisfies the conditions of Theorem \ref{arbitrary}, and if $L$ is defined
by a Galois theoretic property (e.g. $L=K^{ab}$ or $L=K^{solv}$ etc.)
then for any smooth complete curve $X/K$, Corollary \ref{points}
provides a Galois characterization of $L$-rational points.

As a special instance, let us consider the original question of Galois theory
whether any polynomial equation (over ${\bf Q}$) in one variable can be
solved by radicals, and ask the same question for polynomial equations
in two variables, where we now have to assume that the corresponding
curve be geometrically irreducible. The answer to this question is still
unknown (it is equivalent to the question whether ${\bf Q}^{solv}$ is a
PAC-field, cf. [FJ], Problem 10.16(a)). Yet, following the steps of
Galois, we can at least give a group theoretic criterion:
\begin{corollary}
\label{Galois}
Let $X/{\bf Q}$ be a smooth complete curve. Then $X({\bf Q}^{solv})\neq
\emptyset$ iff there is a hensel subgroup $D\leq G_{{\bf Q}(X)}$ with
$G_{{\bf Q}^{solv}}\subseteq pr_{X/{\bf Q}}(D)$.
\end{corollary}
Let us mention that the problem of finding rational points over ${\bf Q}^{solv
}$ is by no means out of date (cf. e.g. Section 4.5 of [SW]).
\subsection{Relating the two fundamental conjectures}
Mochizuki's Theorem A in [Mo] implies both the fundamental conjecture
and the birational fundamental conjecture for smooth hyperbolic
curves over sub-$p$-adic fields (and even for higher-dimensional varieties).
It is, however, not obvious whether, in general, one conjecture implies
the other. The fact that $G_{K(X)}$ can be obtained as inverse limit
of the $\pi_1 (X^\prime)$'s for all Zariski-open $X^\prime\subseteq X$
does not mean that $G_{K(X)}$ has to remember this genealogy.
Yet, for almost all constant fields, it does:
\begin{corollary}
Let $K$ be as in Theorem \ref{arbitrary} and let $X/K$ be a
smooth complete curve over $K$.
Then there is a purely group-theoretic characterization of the
quotients of $G_{K(X)}$ which are $\pi_1(X^\prime )$'s for some
Zarisiki-open $X^\prime\subseteq X$.

In particular, the fundamental conjecture for some
$X^\prime\subseteq X$ over $K$ implies the birational fundamental conjecture
for $X$ over $K$.
\end{corollary}
{\em Proof:}
The previous corollary does not only imply a group-theoretic
characterization of decomposition subgroups $D_P$ of $G_{K(X)}$
w.r.t. points $P\in X(\overline{K})$, but also of the corresponding
inertia subgroups $I_P=D_P\cap\mbox{ker}\,pr_{X/K}$.
The characterization is this: {\em Any subgroup
$I\leq G_{K(X)}$ is an inertia subgroup of $G_{K(X)}$ w.r.t. some
$\overline{K}$-rational point of $X$ iff $I=D\cap\mbox{ker}\,pr_{X/K}$
for some maximal hensel subgroup $D$ of $G_{K(X)}$ containing the
image of a section of $pr_{X/L}:G_{L(X)}\rightarrow G_L$
for some open subgroup $G_L$ of $K$.}

So if ${\cal I}$ is the set of conjugacy classes $[I]$ (in $G_{K(X)}$)
of these group-theoretically described inertia subgroups $I$ of $G_{K(X)}$
one obtains a 1-1 correspondence
$$\begin{array}{rcl}
X(\overline{K}) & \longleftrightarrow & {\cal I}\\
P & \mapsto & [I_P]
\end{array}$$
Note that any two inertia groups belong to the same point iff they are
conjugate.

Now the group-theoretic characterization of fundamental groups is easy:
Let $N\lhd G_{K(X)}$ be a normal subgroup. Then {\em
$G_{K(X)}/N\cong \pi_1(X^\prime)$ for some Zariski-open $X^\prime\subseteq X$
iff $N=\langle I\mid [I]\in {\cal I}\setminus {\cal I}_0\rangle$
for some finite ${\cal I}_0\subseteq {\cal I}$.}

Finally, assume the fundamental conjecture holds for some
Zariski-open $X^\prime\subseteq X$ over $K$,
and let $Y/K$ be a smooth curve over $K$ with $G_{K(Y)}\cong_{G_K}G_{K(X)}$.
We may as well assume that $Y$ is also complete.
Let $N$ be the kernel of $G_{K(X)}\rightarrow \pi_1(X^\prime)$
and let $N^\prime$ be the isomorphic copy of $N$ in $G_{K(Y)}$.
Then $G_{K(Y)}/N^\prime\cong \pi_1(Y^\prime )$
for some Zariski-open $Y^\prime\subseteq Y$ and
$\pi_1(Y^\prime )\cong_{G_K} \pi_1(X^\prime)$.
By assumption this implies $Y^\prime\cong_K X^\prime$ and so $X$ and $Y$
are birationally equivalent over $K$.$\Box$
\subsection{Describing arithmetic properties in Galois-theoretic terms:
the PAC-property, largeness and Hilbertianity}
An immediate consequence of Corollary \ref{points} is the following
group theoretic characterization of PAC-fields and of large fields:
\begin{corollary}
\label{PAC} Let
$K$ be a perfect field satisfying the hypothesis of Theorem \ref{arbitrary}
and let $pr=pr_{{\bf P}^1/K}:G_{K(t)}\rightarrow G_K$
be the canonical projection. Then
\begin{enumerate} \item $K$ is PAC iff every open subgroup
$H$ of $G_{K(t)}$ with $pr(H)=G_K$ contains a hensel subgroup $D$
with $pr (D)=G_K$.
\item
$K$ is large iff every open subgroup $H$ of $G_{K(t)}$ with $pr (H)=G_K$
contains either no or infinitely many pairwise non-conjugate maximal
hensel subgroups $D$ with $pr(D)=G_K$.\end{enumerate}
\end{corollary}
Unlike large fields, Hilbertian fields always satisfy the hypothesis
of Theorem \ref{arbitrary}. Recall that a field $K$ is {\bf separably
Hilbertian} if Hilbert's Irreducibility Theorem holds for separable
polynomials. In characteristic $0$ this is equivalent to Hilbertianity,
and in characteristic $p>0$, $K$ is Hilbertian iff $K$ is imperfect and
separably Hilbertian. If $K$ is separably Hilbertian then so is $K^{perf}$
(cf. [FJ], Section 11.3). It is not known whether the converse holds (cf. [J],
Problem 13). So our Galois theoretic criterion
for separable Hilbertianity can at the moment only be stated for perfect
fields:
\begin{corollary} Let
$K$ be a perfect field satisfying the hypothesis of Theorem \ref{arbitrary}
and let $pr:G_{K(t)}\rightarrow G_K$ be the canonical projection.
Then $K$ is separably Hilbertian iff for every open subgroup $H$ of $G_{K(t)}$
there are infinitely many pairwise non-conjugate maximal hensel subgroups $D$
of $G_{K(t)}$ with $pr(D)=G_K$ such that
$$[D:D\cap H]=[G_{K(t)}:H].$$
\end{corollary}
{\em Proof:}
{\bf `$\Rightarrow$':}
Let $K$ be separably Hilbertian and let $H\leq G_{K(t)}$ be open.
Then the fixed field $F$ of $H$ is of the form $F=K(t,\alpha)$,
where the irreducible polynomial $f(t,Y)$ of $\alpha$ over $K(t)$
is separable and can be chosen in $K[t,Y]$.
Hence there are infinitely many $a\in K$ such that $f(a,Y)\in K[Y]$
is irreducible over $K$. For each such $a$ the $(t-a)$-adic henselisation
$L_a$ of $K(t)$ has the property that
$$[L_a(\alpha ):L_a]=[F:K(t)].$$
$D:=G_{L_a}$ is then a maximal hensel subgroup of $G_{K(t)}$ with $pr(D)=G_K$
such that $D\cap H=G_{L_a(\alpha )}$ and hence
$[D:D\cap H]=[G_{K(t)}:H]$. Finally, distinct $a$'s induce
non-conjugate $D$'s, and the right hand side follows as $K$ is infinite.

{\bf `$\Leftarrow$':}
Assume the right hand side and let $f\in K[t,Y]$ be separable and irreducible.
Choose $\alpha\in\overline{K(t)}$ with $f(t,\alpha)=0$ and let
$H=G_{K(t,\alpha )}$. Choose one of the infinitely many maximal hensel groups
$D$ guaranteed by our assumption avoiding those (up to conjugation) finitely
many corresponding to ramification points of $K(t,\alpha)/K(t)$ or to
zeros of the discriminant of $f$. Then $f(t,Y)$ remains irreducible over
the fixed field of $D$ which, by Corollary \ref{points} is the
henselisation of some $(t-a)$-adic valuation of $K(t)$ with $a\in K$.
By the choice of $D$, $f(a,Y)$ is still separable.
So $f(a,Y)$ is, by Hensel's lemma, irreducible over $K$.$\Box$

Institut f\"{u}r mathematische Logik und Grundlagen der Mathematik\\
Eckerstr. 1, D-79104 Freiburg, Germany\\
e-mail: {\tt Jochen.Koenigsmann@unibas.ch}
\end{document}